\documentclass[
11pt,  reqno]{amsart}
\usepackage{amsmath,amssymb,amscd,amsthm}

\usepackage[hidelinks]{hyperref}

\usepackage{color}

\usepackage{bm}

\allowdisplaybreaks[3]

\newtheorem{theorem}{Theorem} [section]
\newtheorem{lemma}[theorem]{Lemma}
\newtheorem{proposition}[theorem]{Proposition}
\newtheorem{remark}[theorem]{Remark}

\DeclareMathOperator*{\supp}{supp}

\DeclareMathOperator{\Id}{Id}

\newcommand{\Z}{\mathbb{Z}}
\newcommand{\R}{\mathbb{R}}
\newcommand{\C}{\mathbb{C}}
\newcommand{\T}{\mathbb{T}}

\newcommand{\N}{\mathbb{N}}

\renewcommand{\L}{\mathcal{L}}

\renewcommand{\l}{\ell}

\newcommand{\Ml}{\mathcal{M}}

\newcommand{\Hl}{\mathcal{H}}

\newcommand{\cb}{\bm c}
\newcommand{\kb}{\bm k}

\newcommand{\Rf}{\mathfrak{R}}

\newcommand{\fv}{\mathfrak v}
\newcommand{\fw}{\mathfrak w}

\newcommand{\fz}{\mathfrak z}

\newcommand{\al}{\alpha}
\newcommand{\be}{\beta}
\newcommand{\g}{\gamma}

\newcommand{\dl}{\delta}
\newcommand{\Dl}{\Delta}
\newcommand{\eps}{\varepsilon}
\newcommand{\ta}{\theta}

\newcommand{\grad}{\nabla}

\newcommand{\ft}{\widehat}
\newcommand{\Ft}{{\mathcal{F}}}
\newcommand{\wt}{\widetilde}
\newcommand{\cj}{\overline}
\newcommand{\pd}{\partial}
\newcommand{\dx}{\partial_x}

\newcommand{\dt}{\partial_t}

\newcommand{\les}{\lesssim}
\newcommand{\ges}{\gtrsim}

\newcommand{\jb}[1]
{\langle #1 \rangle}

\newcommand{\ind}{\mathbf 1}

\numberwithin{equation}{section}
\numberwithin{theorem}{section}

\begin{document}
\baselineskip = 14pt

\title[Non-existence of solutions to a dNLS system]
{Non-existence of solutions to
a system
\\
of quadratic derivative nonlinear Schr\"odinger equations}

\author[T.~Kondo and M.~Okamoto]
{Toshiki Kondo and Mamoru Okamoto}

\address{
Toshiki Kondo\\
Department of Mathematics\\\
Graduate School of Advanced Science and Engineering\\
Hiroshima University\\
\hfill\break 1-3-1 Kagamiyama, Higashi-Hiroshima, 739-8526 Japan}
\email{d263049@hiroshima-u.ac.jp}

\address{
Mamoru Okamoto\\
Department of Mathematics\\\
Graduate School of Advanced Science and Engineering\\
Hiroshima University\\
\hfill\break 1-3-1 Kagamiyama, Higashi-Hiroshima, 739-8526 Japan}
\email{mokamoto@hiroshima-u.ac.jp}



\date{\today}

\begin{abstract}
We consider the Cauchy problem for a system of quadratic derivative nonlinear Schr\"odinger equations introduced by M.~Colin and T.~Colin (2004) as a model of laser-plasma interaction.
Under a coefficient condition,
we prove that there are initial data admitting no solution in Sobolev spaces,
both on the torus and on the Euclidean space.
\end{abstract}

\maketitle
%
\section{Introduction}
\label{sec:intro}
We consider the Cauchy problem for a system of quadratic derivative nonlinear Schr\"odinger equations:
\begin{equation}\label{eq:cc}
\left\{
\begin{aligned}
&(i\dt+\al \Dl)u=-(\grad \cdot w)v, && (t,x)\in \R \times \Ml^d,\\
&(i\dt+\be \Dl)v=-(\grad \cdot \cj w)u, && (t,x)\in \R \times \Ml^d,\\
&(i\dt+\g \Dl)w=\grad(u \cdot \cj v), && (t,x)\in \R \times \Ml^d,\\
&(u,v,w)|_{t=0} = (u_0,v_0,w_0), && x\in \Ml^d,
\end{aligned}
\right.
\end{equation}
where $\al, \be, \g \in \R \setminus \{0\}$ and $\Ml$ is either $\R$ or $\T :=\R/(2\pi\Z)$.
The unknown functions $u$, $v$, and $w$ are $\C^d$-valued
and
$u_0$, $v_0$, and $w_0$ are initial data.
Here,
$\Dl$ and $\grad$ denote the Laplacian and the gradient on $\Ml^d$, respectively.
This system was introduced by M. Colin and T. Colin \cite{CC} as a model of laser-plasma interaction.

The system \eqref{eq:cc} contains quadratic derivative nonlinearities,
and a derivative loss arising from these terms makes the problem difficult.
For Schr\"odinger equations with derivative nonlinearities,
the well-posedness of the Cauchy problem depends crucially on the structure of the nonlinearity.
Indeed, Christ \cite{Chr} proved that the Cauchy problem for 
\begin{equation*}
(i\dt+\dx^2)u= u \dx u
\end{equation*}
is ill-posed in $H^s(\Ml)$ for any $s\in\R$.
See also \cite{KoOk26} for the real-line case.
In contrast, Gr\"unrock \cite{Gr00} showed that the Cauchy
problem for
\[
\dt u + i \Dl u = \cj{u} (\grad \cdot \cj{u})
\]
is globally well-posed in $H^s(\Ml)$ for $s \ge 0$ and locally well-posed in $H^s(\Ml^d)$ for $d \ge 2$ and $s > \frac d2 -1$.
Related results in the periodic and nonperiodic settings can be found in \cite{Chi02, KO, KoOk26b}.

Next, we recall some known results concerning \eqref{eq:cc}.
The structure of the nonlinearity in \eqref{eq:cc} depends on the coefficients $\al$, $\be$, and $\g$ of the Laplacian.
When
\[
(\al - \g)(\be + \g) \neq 0,
\]
local well-posedness in
\[
\Hl^s (\Ml^d) := (H^s(\Ml^d))^d \times (H^s(\Ml^d))^d \times (H^s(\Ml^d))^d
\]
was established on both $\R^d$ and $\T^d$;
see Tables 1--3 in \cite{HKO25}.
These results were obtained by an iteration argument in Fourier restriction norm spaces.

On the other hand,
when
\[
(\al - \g)(\be + \g) = 0,
\]
Hirayama \cite{Hir14} proved that the flow map is not $C^2$ in $\Hl^s(\R)$ for any $s\in\R$.
This shows that well-posedness does not follow from an iteration argument.
When $\be +\g \neq 0$, however,
Hirayama, Kinoshita, and the second author \cite{HKO22} showed local well-posedness in $\Hl^s(\R^d)$ for $s>\frac d2+3$ by a modified energy argument.
Their argument can also be applied to the periodic case.
See also \cite{Ak25} for the case $\al - \g = 0$ and $\be + \g \neq 0$.

In contrast, the case $\be +\g = 0$ exhibits different behavior.
Indeed, in the periodic setting, 
Hirayama, Kinoshita, and the second author \cite{HKO25} proved that the flow map of \eqref{eq:cc} fails to be continuous in $\Hl^s(\T^d)$ for any $s\in \R$.

The purpose of this paper is to show that, for both $\Ml = \R$ and $\Ml = \T$,
the condition 
\[ 
\be + \g = 0 
\] 
leads to the non-existence of solutions to \eqref{eq:cc} in $\Hl^s(\Ml^d)$ for $s > \frac d2 + 1$.
The regularity assumption $s> \frac d2+1$ allows us to control the nonlinear remainders using the algebra property of $H^{s-1}(\Ml^d)$.
Our main result is as follows.

\begin{theorem}
\label{thm:nonexist}
Assume that
\begin{equation}
\be+\g=0.
\label{cond0}
\end{equation}
Let $\Ml$ be $\R$ or $\T$, and let $s>\frac d2+ 1$.
For every $\eps>0$, there exist initial data
$(u_0,v_0,w_0)\in \Hl^s(\Ml^d)$ such that
\[
\|u_0\|_{(H^s)^d}+\|v_0\|_{(H^s)^d}+\|w_0\|_{(H^s)^d}<\eps
\]
and \eqref{eq:cc} on $\Ml^d$ admits no distributional solution in
$C([0,T];\Hl^s(\Ml^d))$ for any $T>0$.
Moreover, the data can be chosen with
\[
u_0\in\bigcap_{m=0}^{\infty} (H^m(\Ml^d))^d,
\quad
v_0=0.
\]
\end{theorem}

Under \eqref{cond0},
the linear parts of the equations for $v$ and $\cj w$ coincide.
For the constant background $u= (c,0, \dots, 0)$, the corresponding Fourier system has a growing component with rate $|c||k|$ for frequencies along the $x_1$-axis;
see \eqref{eq:matrix}--\eqref{eq:growth}.
This yields a necessary regularity condition for the initial data, in the spirit of the scalar argument in \cite{KO}.
While the result of \cite{HKO25} established norm inflation on $\T^d$,
our result demonstrates a different ill-posedness phenomenon:
the non-existence of local solutions for arbitrarily small initial data on both $\T^d$ and $\R^d$.

We now outline the proof.
Set
\[
q := \cj w.
\]
Taking the complex conjugate of the third equation in \eqref{eq:cc} and using \eqref{cond0},
we obtain
\[
(i\dt+\be \Dl)q=- \grad (\cj u \cdot v).
\]
Thus, $v$ and $q$ have the same linear part
and satisfy
\begin{equation}\label{eq:split}
\begin{aligned}
\dt v&=i\be \Dl v+ i (\grad \cdot q) u_0 + R^{(1)},\\
\dt q&=i\be \Dl q+ i \grad (\cj {u_0} \cdot v) + R^{(2)},
\end{aligned}
\end{equation}
where
\begin{equation}\label{eq:R12}
R^{(1)} := i (\grad \cdot q)(u-u_0),
\quad
R^{(2)} := i \grad( (\cj{u-u_0}) \cdot v).
\end{equation}
Formally neglecting the dispersive terms and the remainders $R^{(1)}$ and $R^{(2)}$,
we arrive at the leading-order system
\[
\begin{aligned}
\dt v&= i u_0 (\grad \cdot q),\\
\dt q&= i \grad (\cj{u_0} \cdot v).
\end{aligned}
\]
When $\Ml = \T$ and $u_0$ is a constant vector $\cb \in \C^d$,
this system becomes
\[
\begin{aligned}
\dt v&= i \cb (\grad \cdot q),\\
\dt q&= i \grad (\cj{\cb} \cdot v).
\end{aligned}
\]
Taking $\cb = (c, 0, \cdots, 0)$ and writing $v_1$ and $q_\l$ for the first component of $v$ and the $\l$-th component of $q$, respectively,
we obtain
\[
\begin{aligned}
\dt \ft v_1 (t, \kb) &= - c \sum_{\l = 1}^d k_\l \ft q_\l (t,\kb),\\
\dt \ft q_\l (t,\kb) &= - \cj  c k_\l \ft v_1 (t,\kb)
\end{aligned}
\]
for $\kb = (k_1, \dots, k_d) \in \Z^d$.
In particular, for $\kb = (k_1, 0, \cdots, 0)$ with $k_1 \in \Z \setminus \{0\}$, we have
\[
\begin{aligned}
\dt^2 \ft q_1 (t,\kb) - |c|^2 |k_1|^2 \ft q_1(t,\kb) = 0.
\end{aligned}
\]
The operator $\dt^2 - |c|^2 |k_1|^2$ is the Fourier representation of the elliptic operator $\dt^2 + |c|^2 \pd_{x_1}^2$. 
This suggests that the existence of a solution with $v_0=0$ imposes additional Sobolev regularity on the initial data in the $x_1$-direction.

These observations suggest that the leading-order dynamics of \eqref{eq:cc} under \eqref{cond0} exhibit elliptic-type behavior.
To turn this heuristic into a rigorous non-existence argument,
we must control the nonlinear terms.

For scalar semilinear Schr\"odinger equations,
the authors \cite{KO} proved the non-existence of solutions using a gauge transformation.
In the present system,
an analogous transformation would be matrix-valued and would introduce additional commutator terms,
so the argument in \cite{KO} does not apply directly.

Instead, we exploit the factor of $t$ gained by estimating the difference between the solution and its initial data.
This allows us to control the nonlinear remainders on a sufficiently short time interval.
We then exploit the elliptic-type behavior to show that the existence of a solution to \eqref{eq:cc}  under \eqref{cond0} forces additional regularity of the initial data in the $x_1$-direction.
Choosing initial data that fail to have this additional regularity yields the desired non-existence result.

When $\Ml = \R$,
we cannot choose a nonzero constant vector $u_0 = \cb \in \C^d$,
since such a function does not belong to $(H^s(\R^d))^d$.
Instead,
we use the wave-packet method to extract the worst interaction.
Wave packets have been used to establish necessary conditions for the existence of solutions to linear Schr\"odinger equations \cite{Mi85};
see also \cite{JeOh26} for a quasilinear case.
In \cite{KoOk26b}, the authors used this method to prove the non-existence of solutions to a semilinear  Schr\"odinger equation in $H^s(\R)$.
We adapt the method to the system \eqref{eq:cc}.
On a region where $u_0$ is a real constant vector $\cb$,
set
\[
z := {v_1 -q_1}, \quad R := {R_{1}^{(1)} - R_1^{(2)}}.
\]
Then,
\eqref{eq:split} yields that
\[
\dt z = i \be \Dl z - i u_{0,1} \pd_{x_1} z + i u_{0,1} \sum_{\l = 2}^d \pd_{x_\l} q_\l + R
\]
on that domain.
Thus,
we obtain a scalar Schr\"odinger equation with an additional coupling term,
which can be controlled as an error.
This allows us to adapt the method in \cite{KoOk26b}.
Unlike the single wave-packet construction used in \cite{KoOk26b}, 
we employ infinitely many wave packets to show that the existence of a solution forces the initial data to have higher regularity.
See Remark \ref{rem:reg} for further details.

\subsection{Notation}
We conclude this section by fixing notation.
We denote by $\N$ the set of positive integers and by $\N_0$ the set of non-negative integers.
For $x \in \C^d$, we write $x_j$ for its $j$-th component.

On $\T^d$, we use normalized Lebesgue measure and define the Fourier coefficients by
\begin{align*}
\int_{\T^d} f dx &:=\frac1{(2\pi)^d}\int_{[- \pi, \pi]^d} f(x) dx,
\\
\ft f(\kb) = \Ft [f] (\kb) &:=\int_{\T^d} f(x)e^{-i \kb \cdot x} dx
\end{align*}
for $\kb \in \Z^d$.
On $\R^d$, we use Lebesgue measure and the unitary Fourier transform
\[
\ft f(\xi) = \Ft [f] (\xi) := \frac 1{(2\pi)^{\frac d2}} \int_{\R^d} f(x)e^{-ix \cdot \xi} dx
\]
for $\xi \in \R^d$.

For $r \in \R$,
let $\jb{\dx}^r$ denote the Fourier multiplier with the symbol $\jb{\xi}^r$
($\xi \in \Z^d$ on $\T^d$),
where
$\jb{\xi}:=\sqrt{1+|\xi|^2}$.
The Sobolev norm is given by
\[
\|f\|_{H^r}=\|\jb{\dx}^{r}f\|_{L^2}.
\]
For a $\C^d$-valued function $f$,
we write  $f_j$ for its $j$-th component and define its Sobolev norm by 
\[
\|f\|_{(H^r)^d} = \sum_{j=1}^d \| f_j \|_{H^r}.
\]
We abbreviate
$L_T^\infty (H^r)^d := L^\infty([0,T]; (H^r(\Ml^d))^d)$.

We write $A \les B$ if $A\le CB$ for some constant $C>0$,
and $A\sim B$ if $A\les B$ and $A\ges B$.

\section {Non-existence}

In this section,
we prove Theorem \ref{thm:nonexist}.
Let
$s>\frac d2 + 1$ and let
\[
\begin{aligned}
u_0 &= (u_{0,1}, 0, \dots, 0) \in ( H^{s+1}(\Ml^d) )^d, \quad
v_0 = 0, \\
w_0 &= (w_{0,1}, 0, \dots, 0) \in ( H^s(\Ml^d) )^d.
\end{aligned}
\]
Suppose that there exists a solution
\[
(u,v,w) \in C([0,T]; \Hl^s(\Ml^d))
\]
to \eqref{eq:cc} with \eqref{cond0} for some $T>0$.
Restricting the interval, we may assume $T \le 1$.
Set
\[
M = 1 + \|u_0\|_{(H^{s+1})^d} + \|u \|_{L_T^\infty (H^s)^d} + \|v \|_{L_T^\infty (H^s)^d} + \|w \|_{L_T^\infty (H^s)^d} < \infty.
\]
All constants in the rest of the proof may depend on $M$, $s$, and the fixed coefficients.

First, we estimate the remainder terms $R^{(1)}$ and $R^{(2)}$ defined in \eqref{eq:R12}.
Since $s-1>\frac d2$,
the space $H^{s-1}(\Ml^d)$ is a Banach algebra:
\[
\|fg\|_{H^{s-1}}\les\|f\|_{H^{s-1}}\|g\|_{H^{s-1}}.
\]
The nonlinearities in \eqref{eq:cc} are continuous as $( H^{s-1}(\Ml^d))^d$-valued functions of time.
Hence,
the mild form of the $v$ equation, using $v|_{t=0}=0$, is
\[
v(t)= i \int_0^t e^{i(t-\tau)\be\Dl}
( u (\tau) (\grad \cdot q(\tau) ) ) d\tau.
\]
By the algebra estimate, we have
\begin{equation}
\|v(t)\|_{( H^{s-1} )^d}
\les\int_0^t\|u (\tau)\|_{(H^{s-1})^d}\|q(\tau)\|_{(H^s)^d} d\tau
\les t.
\label{eq:vt1}
\end{equation}
Using $u|_{t=0}=u_0$, the mild form of the $u$ equation is 
\[
\begin{aligned}
u(t) 
&= e^{it \al\Dl} u_0 + i \int_0^t e^{i(t-\tau) \al \Dl}
\big( (\grad \cdot w(\tau) ) v(\tau) \big) d\tau.
\end{aligned}
\]
By the algebra estimate, we have
\begin{equation}
\begin{aligned}
\| u(t) - u_0 \|_{(H^{s-1})^d}
&\le
\| ( e^{it \al\Dl} - \Id ) u_0 \|_{(H^{s-1})^d} 
\\
&\quad+
\int_0^t \|q(\tau)\|_{(H^{s})^d}\|v(\tau)\|_{(H^{s-1})^d} d\tau
\\
&
\les 
t \| u_0 \|_{(H^{s+1})^d} + t^2
\les t.
\end{aligned}
\label{eq:ut1}
\end{equation}
From \eqref{eq:R12}, \eqref{eq:ut1}, and \eqref{eq:vt1},
we have
\begin{equation}\label{eq:e-orders}
\begin{aligned}
\|R^{(1)}(t)\|_{(H^{s-1})^d}
&\les
\| u(t) - u_0 \|_{(H^{s-1})^d} \| q(t) \|_{(H^{s})^d}
\les 
t, \\
\| R^{(2)}(t) \|_{(H^{s-2})^d}
&\les
\| u(t) - u_0 \|_{(H^{s-1})^d} \| v(t) \|_{(H^{s-1})^d}
\les
t^2.
\end{aligned}
\end{equation}

By considering each component of $v$ and $q$ separately, \eqref{eq:split} becomes
\begin{equation}
\begin{aligned}
\dt v_1 &= i \be \Dl v_1 + i u_{0,1} (\grad \cdot q) + R^{(1)}_{1},
\\
\dt q_\l &= i \be \Dl q_\l + i \pd_{x_\l} (\cj {u_{0,1}} v_1) + R^{(2)}_{\l}
\end{aligned}
\label{sys1}
\end{equation}
for $1 \le \l \le d$.
Since $w_0 = (w_{0,1}, 0, \dots, 0)$, the mild form of the equations for $q_\l$ is
\begin{equation*}
\begin{aligned}
q_\l(t) = \int_0^t e^{i(t-\tau)\be \Dl} ( i \pd_{x_\l} (\cj {u_{0,1}} v_1(\tau)) + R^{(2)}_{\l}(\tau) ) d \tau 
\end{aligned}
\end{equation*}
for $2 \le \l \le d$.
From the algebra estimate, \eqref{eq:vt1}, and \eqref{eq:e-orders},
we obtain
\begin{equation}
\begin{aligned}
\| q_\l (t) \|_{H^{s-2}}
&\les
\int_0^t \big( \| \cj {u_{0,1}}\|_{H^{s-1}} \| v_1(\tau) \|_{H^{s-1}} +  \| R^{(2)}_{\l}(\tau) \|_{H^{s-2}} \big) d \tau 
\\
&\les
t^2 + t^3 \les t^2
\end{aligned}
\label{sys3}
\end{equation}
for every $2 \le \l \le d$.

\subsection{On the torus}
Let
\[
E := \{ \kb \in \Z^d \mid \kb = a (1, 0, \cdots, 0) \text{ for some } a \in \Z \setminus\{0\} \},
\]
and define $\ind_E$ to be the characteristic function of $E$:
\[
\ind_E(\kb) := 
\begin{cases}
1, & \text{if } \kb \in E, \\
0, & \text{otherwise}.
\end{cases}
\]
Then, we define $\Pi_1$ as
\[
\Pi_1 = \Ft^{-1} \ind_E(\kb) \Ft. 
\]

On $\T^d$, we prove the following.

\begin{proposition}
\label{prop:necessary}
Let $\Ml = \T$, $s>\frac d2 + 1$, $T>0$, and $c\in \C \setminus\{0\}$.
If \eqref{eq:cc} with \eqref{cond0} has a solution $(u,v,w) \in C([0,T];\Hl^s(\T^d))$ with
\begin{equation}
\begin{aligned}
u_0 &= \cb := (c, 0, \cdots, 0), 
\quad
v_0 = 0,
\quad
w_0 = (w_{0,1}, 0, \cdots, 0),
\end{aligned}
\label{iniT}
\end{equation}
then
\[
\Pi_1 w_{0,1} \in H^{s+\dl}(\T^d)
\]
for every $0<\dl<1$.
\end{proposition}

Once we obtain Proposition \ref{prop:necessary},
Theorem \ref{thm:nonexist} on $\T^d$ follows from the fact that there exists $w_{0,1} \in H^s(\T^d)$ such that  
$\Pi_1 w_{0,1} \notin H^{s+\dl}(\T^d)$ for $0<\dl<1$.
Indeed,
for small $\eps>0$,
we take
\[
u_0 (x) = \Big(\frac \eps2, 0, \cdots, 0 \Big), \quad
w_0(x)= \Big( \eps^2 \sum_{j\in\N}\frac{2^{-js}}{j}e^{-i 2^j x_1}, 0, \cdots, 0 \Big).
\]
Then,
$w_{0,1}$ satisfies
\[
\| \Pi_1 w_{0,1} \|_{H^s}^2
\sim
\eps^4 \sum_{j\in\N}\frac1{j^2}
< \frac{\eps^2}4.
\]
Moreover,
\begin{equation*}
\|\Pi_1 w_{0,1}\|_{H^{s+\dl}}^2
\sim
\eps^4
\sum_{j\in\N}\frac{2^{2j\dl}}{j^2}=\infty
\end{equation*}
for $\dl>0$.
This contradicts Proposition \ref{prop:necessary},
which proves Theorem \ref{thm:nonexist} on $\T^d$.

\begin{proof}[Proof of Proposition \ref{prop:necessary}.]
Define
\begin{align*}
\fv(t) &:=e^{-it\be\Dl}v(t), \\
\fw(t) &:=e^{-it\be\Dl}q(t), \\
\wt R^{(j)}(t) &:=e^{-it\be\Dl}R^{(j)}(t)
\end{align*}
for $j=1,2$.
Then,
\eqref{sys1} with the initial data \eqref{iniT} becomes
\begin{equation}
\begin{aligned}
\dt {\ft \fv_1 (t, \kb)} 
&= - c k_1 \ft \fw_1 (t, \kb) - c \sum_{\l=2}^d k_\l \ft \fw_\l (t, \kb)  + \Ft [\wt R^{(1)}_{1}] (t, \kb) , \\
\dt {\ft \fw_1 (t, \kb)} 
&= - \cj c k_1 \ft \fv_1 (t, \kb)  + \Ft [\wt R^{(2)}_{1}] (t, \kb).
\end{aligned}
\label{eq:matrix}
\end{equation}
Define $\fz$ by
\[
\ft \fz(t,\kb) := \frac{\cj c}{|c|} \ft \fv_1(t,\kb) - \frac {k_1}{|k_1|} \ft \fw_1(t,\kb)
\]
for $\kb \in E$,
and $\ft \fz(t,\kb)=0$ for $\kb \notin E$.
Then, $\fz\in C([0,T];H^s(\T^d))$.
From \eqref{eq:matrix},
we have
\begin{equation}
\dt\ft \fz(t, \kb)
= |c| |k_1|  \ft \fz(t,\kb) 
+ \ft \Rf(t,\kb)
\label{eq:growth}
\end{equation}
for $\kb \in E$,
where
\[
\ft \Rf(t,\kb)
=
\frac{\cj c}{|{c}|} \Ft [\wt R^{(1)}_{1}](t,\kb) - \frac {k_1}{|k_1|} \Ft [\wt R^{(2)}_{1}](t,\kb).
\]
We set $\ft \Rf(t, \kb) = 0$ when $\kb \notin E$.
Hence,
we have
\[
\dt (\ind_E(\kb) \ft \fz(t,\kb)) = |c||\kb| \ind_E (\kb) \ft \fz(t, \kb) + \ind_E(\kb) \ft \Rf(t, \kb).
\]

Duhamel's formula  gives
\[
\ind_E(\kb) \ft \fz(T,\kb)
=
e^{|c| |\kb| T} \ind_E(\kb) \ft \fz(0,\kb)
+ \int_0^T
e^{|c| |\kb| (T-t)} \ind_E(\kb) \ft \Rf(t,\kb) dt.
\]
Hence, we have
\begin{equation}
\begin{aligned}
&\ind_E (\kb) \ft \fz(0,\kb)
\\
&= e^{- |c| |\kb| T} \ind_E(\kb) \ft \fz(T,\kb)
-
\int_0^T e^{- |c| |\kb|t} \ind_E (\kb) \ft \Rf(t,\kb) dt.
\end{aligned}
\label{Duhfz}
\end{equation}
The first term on the right-hand side decays exponentially,
because $\fz(T)\in H^s(\T^d)$ and $T>0$.
Here and below, $|\grad|$ denotes the Fourier multiplier with symbol $|\kb|$.
For $r'\ge0$, $r\in\R$, and $0<t\le1$, we use
\[
\|e^{-|c|t |\grad|} \Pi_1 f\|_{H^{r+r'}}
\les t^{-r'}\|\Pi_1 f\|_{H^r},
\]
which follows from $\sup_{ \kb \in E}\jb{\kb}^{r'}e^{-|\cb|t|\kb|}\les t^{-r'}$.

From \eqref{eq:e-orders},
we obtain
\begin{equation}
\begin{aligned}
\|e^{-|c| t |\grad|} \Pi_1 \wt R^{(1)}(t)\|_{(H^{s+\dl})^d}
&\les t^{-(1+\dl)}\|R^{(1)}(t)\|_{(H^{s-1})^d}
\les t^{-\dl},
\\
\|e^{-|c| t |\grad|} \Pi_1 \wt R^{(2)}(t)\|_{(H^{s+\dl})^d}
&\les t^{-(2+\dl)}\|R^{(2)}(t)\|_{(H^{s-2})^d}
\les t^{-\dl}
\end{aligned}
\label{bddT1}
\end{equation}
for $\dl>0$.
When $0 < \dl<1$,
it follows from \eqref{Duhfz}, \eqref{bddT1}, and 
\[
\int_0^T t^{-\dl}dt=\frac{T^{1-\dl}}{1-\dl}<\infty
\]
that
$\Pi_1 \fz(0)\in H^{s+\dl}(\T^d)$.
Note that
\[
\ind_E (\kb) \ft \fz(0,\kb)
=
- \frac {k_1}{|k_1|} \ind_E (\kb) \Ft [\cj {w_{0,1}}](\kb)
=
- \frac {k_1}{|k_1|} \ind_E (\kb) \cj{\Ft [w_{0,1}](-\kb)},
\]
because $v|_{t=0}=0$.
We obtain
\[
\Pi_1 w_{0,1} \in H^{s+\dl}(\T^d),
\]
which completes the proof of Proposition \ref{prop:necessary}.
\end{proof}

\subsection{On the Euclidean space}

The periodic argument uses a nonzero constant function as the initial datum for $u$,
but such a function does not belong to $(H^s(\R^d))^d$.
We replace that constant by infinitely many separated plateaus whose heights tend to zero and whose lengths in the $x_1$-direction tend to infinity.

Set
\[
(f,g)_{L^2}=\int_{\R^d} f\cj g dx.
\]
We use the same notation for the corresponding Sobolev duality pairing whenever needed.

\begin{lemma}\label{lem:mod}
Let $\chi\in C_c^\infty(\R^d)$, $m \in \N_0$, $r\in\R$, and $\nu \in \N_0^d$ satisfy 
\[
\nu_1 + \cdots + \nu_d = m.
\]
For $A, a \in\R^d$, one has
\[
\|e^{iA \cdot x} \dx^\nu \chi(x-a) \|_{H^r}
\les \jb{A}^r,
\]
where $\dx^\nu = \pd_{x_1}^{\nu_1} \cdots \pd_{x_d}^{\nu_d}$.
The constant may depend on $r$, $m$, and $\chi$, but is independent of $A$ and $a$.
\end{lemma}

\begin{proof}
Translation changes only the phase of the Fourier transform.
Peetre's inequality gives
\[
\jb{A+\xi}^{r}
\les
\jb{A}^{r}\jb{\xi}^{|r|}
\]
for $\xi \in \R^d$.
Apply this inequality to the Fourier transform of $\dx^\nu \chi$,
which is rapidly decreasing.
\end{proof}

We construct the initial data.
We fix real-valued functions $\theta,\chi\in C_c^\infty(\R^d)$ with
\begin{gather*}
0 \le \theta \le 1,
\qquad
\theta = 1\ \text{on} \; [-2,2]^d,
\qquad
\supp \theta \subset (-3,3)^d,
\\
\supp \chi \subset (-1, 1)^d,
\qquad
\|\chi\|_{L^2}=1.
\end{gather*}
For $j\in\N$, set
\begin{equation}
\begin{aligned}
c_j &:= 2^{- 3 j},
\quad
x_\ast^{(j)} := 10 \sum_{k=1}^j 2^{4 k},
\quad
N_j := 2^{11 j},
\end{aligned}
\label{eq:constants}
\end{equation}
Note that 
the $d$-dimensional rectangles
\[
( x_\ast^{(j)}- 3 \cdot 2^{4j} , x_\ast^{(j)}+3 \cdot 2^{4j} ) \times (-3,3)^{d-1}
\]
are pairwise disjoint.

Let $\{ a_j \}_{j \in \N}$ be a summable sequence of positive numbers.
For an integer $J\ge1$ (to be fixed later),
define
\begin{equation}
\begin{aligned}
u_{0,1} (x) &:= \sum_{j\ge J} c_j \ta \Big( \frac{x_1- x_\ast^{(j)}}{2^{4j}}, x_2, \cdots ,x_d \Big), 
\\
w_{0,1}(x) &:= \sum_{j\ge J} a_j N_j^{-s} e^{-i N_j (x_1-x_\ast^{(j)})} \chi(x_1 - x_\ast^{(j)}, x_2, \cdots, x_d), 
\\
u_0(x) &= (u_{0,1}(x), 0, \cdots, 0 ),
\\
v_0(x) &= 0,
\\
w_0(x) &= (w_{0,1}(x), 0, \cdots, 0).
\end{aligned}
\label{eq:ini}
\end{equation}
For each nonnegative integer $m$,
we have
\[
\| u_{0,1} \|_{H^m}
\les
\sum_{j \ge J}
c_j \cdot 2^{2 j} =
\sum_{j \ge J}
2^{- j}
< \infty.
\]
Thus,
$u_0 \in\bigcap_{m\ge0} (H^m (\R^d))^d$,
and $u_0= (c_j, 0, \cdots, 0)$ on
the $d$-dimensional rectangle
\begin{equation}
X_j := ( x_\ast^{(j)}- 2 \cdot 2^{4j} , x_\ast^{(j)}+ 2 \cdot 2^{4j} ) \times (-2,2)^{d-1}.
\label{IX}
\end{equation}
Lemma~\ref{lem:mod} with $A = -N_j (1,0,\cdots,0)$ gives
\[
\|w_{0,1}\|_{H^s} \les \sum_{j \ge J} a_j < \infty.
\]
In particular, $(u_0, 0, w_0)\in \Hl^s(\R^d)$.

Taking an integer $m\ge s$ in the estimate for $u_0$, we obtain a constant $C_s>0$, independent of $J$, such that
\begin{equation*}
\|u_0\|_{(H^s)^d}+\|w_0\|_{(H^s)^d}
\le C_s\sum_{j \ge J} (2^{- j}+a_j) <\eps,
\end{equation*}
provided $J$ is sufficiently large. Fix such a $J$ and the corresponding initial data.

Note that $u_0$ defined in \eqref{eq:ini} is not integrable.
Indeed,
\[
\int_{\R^d} |u_0| dx
= \int_{\R^d} \theta dx \sum_{j \ge J} c_j 2^{4j}
= \int_{\R^d} \theta dx \sum_{j \ge J} 2^{j}
= \infty.
\]

\begin{proposition}
\label{prop:necessaryR}
Let $\Ml=\R$, $s>\frac d2 + 1$, and $T>0$.
If \eqref{eq:cc} with \eqref{cond0} and \eqref{eq:ini} has a solution $(u,v,w) \in C([0,T];\Hl^s(\R^d))$,
then
\begin{equation}\label{eq:e-key}
a_j \les
e^{-j}
\end{equation}
for large $j \in \N$.
\end{proposition}

Theorem \ref{thm:nonexist} on $\R^d$
follows from Proposition \ref{prop:necessaryR} with $a_j = 2^{-j}$.

\begin{remark}
\label{rem:reg}
\rm
It follows from Lemma \ref{lem:mod}, \eqref{eq:ini}, \eqref{eq:constants}, and \eqref{eq:e-key} that
\[
\| w_{0,1} \|_{H^{s+\dl}}
\les
\sum_{j \ge J} a_j N_j^{\dl}
\les
\sum_{j \ge J} e^{-j} 2^{11 \dl j}
< \infty
\]
provided that $\dl \le \frac 1{11}$.
In particular,
we have $w_{0,1} \in H^{s+\dl}(\R^d)$.
Since Theorem \ref{thm:nonexist} follows directly from Proposition \ref{prop:necessaryR},
we do not pursue the optimal regularity gain here.
\end{remark}

\begin{proof}[Proof of Proposition \ref{prop:necessaryR}]
Recall that $X_j$ is defined in \eqref{IX}.
Since 
$u_0(x)= (c_j, 0, \cdots, 0)$ on $X_j$,
the system \eqref{sys1} becomes
\begin{align*}
\dt v_1&=i \be \Dl v_1 + i c_j (\grad \cdot q) + R_{1}^{(1)},\\
\dt q_\l&=i\be \Dl q_\l + i {c_j} \pd_{x_\l} v_1 + R_{\l}^{(2)}
\end{align*}
for $1 \le \l \le d$.
Set
\[
z := v_1 - q_1,
\quad
R := R_{1}^{(1)} - R_{1}^{(2)},
\]
and
\[
\L_j := i\be \Dl - i c_j \pd_{x_1}.
\]
On $X_j$ as in \eqref{IX},
the following identity holds in the sense of distributions:
\begin{equation}\label{eq:e-local}
\dt z=\L_{j} z + i c_j \sum_{\l=2}^d \pd_{x_\l} q_\l + R.
\end{equation}

Define
\begin{equation}
\begin{aligned}
\wt \chi_j(t,x) &:= \chi (x_1 - x_\ast^{(j)} - 2 \be N_j t, x_2, \cdots, x_d),
\\
p_j (t,x) &:= \exp \big( - c_j N_j t - i \be N_j^2 t + i N_j (x_1 - x_\ast^{(j)}) \big),
\\
h_j(t,x) &:=
p_j (t,x)
\wt \chi_j (t,x).
\label{eq:e-packet}
\end{aligned}
\end{equation}
We also set
\begin{equation}\label{eq:e-time}
t_j
:=\frac{j}{c_j N_j}
= j \cdot 2^{- 8 j}.
\end{equation}
One has $t_j\le T$ for all sufficiently large $j \in \N$.

It follows from \eqref{eq:constants} that
\[
2|\be|N_j t_j
\le
2 |\be| j \cdot 2^{3j}
\le \frac{2^{4 j}}2
\]
for large $j \in \N$.
From \eqref{eq:e-packet}, \eqref{eq:e-time}, and \eqref{IX},
we have
\begin{equation}
\supp
h_j(t, \cdot)
\subset
X_j
\label{eq:e-supph}
\end{equation}
for large $j \in \N$ and $0\le t\le t_j$.

The adjoint operator of $\L_j$ is
\begin{align*}
\L_j^* &:= -i \be \Dl - i c_j \pd_{x_1}.
\end{align*}
By \eqref{eq:e-packet},
a direct calculation gives 
\[
(\dt + \L_j^\ast) p_j (t,x) = 0.
\]
From \eqref{eq:e-packet},
we have
\begin{equation}
\begin{aligned}
E_j(t,x)
&:= (\dt+\L_{j}^\ast) h_j (t,x)
\\
&=
- i p_j (t,x) ( \be (\Dl \wt \chi_j) (t,x) +  c_j (\pd_{x_1} \wt \chi_j) (t,x) ) .
\end{aligned}
\label{eq:e-resid}
\end{equation}
It follows from  Lemma \ref{lem:mod} that
\begin{equation}\label{eq:e-norms}
\|h_j(t)\|_{H^r}
+
\|E_j(t)\|_{H^r}
\les N_j^r e^{-c_j N_j t}
\end{equation}
and
\begin{equation}
\|\pd_{x_\l} h_j(t)\|_{H^r} \les N_j^r e^{-c_j N_j t}
\label{hnorm}
\end{equation}
for $r \in \R$ and $2 \le \l \le d$.

The distributional equations imply
\[
u,v,q\in C^1 \big( [0,T]; (H^{s-2}(\R^d))^d \big).
\]
From \eqref{eq:e-local}, \eqref{eq:e-supph}, and \eqref{eq:e-resid},
we  obtain
\begin{equation}\label{eq:e-pair}
\begin{aligned}
&(z(t_j),h_j(t_j))_{L^2}-(z(0),h_j(0))_{L^2}\\
&
= \int_0^{t_j} \Big( i c_j \sum_{\l = 2}^d (\pd_{x_\l} q_\l (t), h_j(t))_{L^2} 
\\
&\hspace*{50pt}
+ (R(t),h_j(t))_{L^2} + (z(t),E_j(t))_{L^2} \Big) dt.
\end{aligned}
\end{equation}
Since
\[
z(0,x)= - \cj {w_{0,1}(x)} = - \sum_{j\ge J} a_j N_j^{-s} e^{i N_j (x_1-x_\ast^{(j)})} \chi(x_1 - x_\ast^{(j)}, x_2, \cdots, x_d),
\]
it follows from \eqref{eq:ini}, \eqref{eq:e-packet}, and \eqref{eq:e-supph} that
\begin{equation}\label{eq:e-init}
| (z(0),h_j(0))_{L^2} |
={a_j N_j^{-s}}.
\end{equation}

From $z \in C([0,T]; H^s(\R^d))$ and \eqref{eq:e-norms},
we have
\begin{align*}
N_j^s|(z(t_j),h_j(t_j))_{L^2}|
&\le
N_j^s \| z(t_j) \|_{H^s} \| h_j(t_j) \|_{H^{-s}}
\les e^{- c_j N_jt_j},
\\
N_j^s\int_0^{t_j}|(z(t),E_j(t))_{L^2}| dt
&\les
\int_0^{t_j} e^{-c_j N_j t} dt
\les
\frac 1{ c_j N_j}.
\end{align*}
Using \eqref{eq:e-orders} and \eqref{eq:e-norms} gives
\begin{align*}
N_j^s\int_0^{t_j}|(R^{(1)}_{1}(t), h_j(t))_{L^2}|dt
&\les
N_j \int_0^{t_j}
t e^{- c_j N_jt} dt
\les
\frac{1}{c_j^2 N_j},\\
N_j^s \int_0^{t_j}|( R^{(2)}_{1}(t) , h_j(t) )_{L^2} |dt
&\les
N_j^2\int_0^{t_j} t^2 e^{-c_j N_jt}dt
\les
\frac{1}{c_j^3 N_j}.
\end{align*}
For $2 \le \l \le d$, integration by parts, \eqref{sys3}, and \eqref{hnorm} yield that
\begin{equation*}
\begin{aligned}
c_j N_j^s \int_0^{t_j} |(\pd_{x_\l} q_\l(t), h_j(t))_{L^2}| dt 
&=
c_j N_j^s \int_0^{t_j} |(q_\l(t), \pd_{x_\l} h_j(t))_{L^2}| dt
\\
&\les
c_j N_j^2 \int_0^{t_j} t^2 e^{-c_j N_j t} dt
\\
&\les
\frac 1{c_j^2 N_j}.
\end{aligned}
\end{equation*}

Combining these estimates with \eqref{eq:e-pair}, \eqref{eq:e-init}, \eqref{eq:constants}, and \eqref{eq:e-time},
we have
\[
a_j \les e^{-c_j N_j t_j}+\frac{1}{c_j^3 N_j}
=
e^{-j}
+2^{-2 j}
\les
e^{-j},
\]
which shows \eqref{eq:e-key}.
\end{proof}

\mbox{}

\noindent
{\bf 
Acknowledgements.}
This work was
supported by JSPS KAKENHI Grant Number
JP26K00613 and JST SPRING Grant Number JPMJSP2132.


\begin{thebibliography}{99}
\bibitem{Ak25}
K. Akase, 
{\it Local well-posedness for a system of quadratic derivative nonlinear Schr\"odinger equations}, 
SUT J. Math. {\bf 62} (2026), no.~1, 63--110.

\bibitem{Chi02}
H. Chihara, 
{\it The initial value problem for Schr\"odinger equations on the torus}, 
Int. Math. Res. Not. {\bf 2002}, no.~15, 789--820.

\bibitem{CC}
M.~Colin, T.~Colin,
\textit{On a quasilinear Zakharov system describing laser-plasma interactions},
Differential Integral Equations \textbf{17} (2004), no. 3-4, 297--330.

\bibitem{Chr}
M. Christ,
{\it Illposedness of a Schr\"odinger equation with derivative nonlinearity},
preprint
(\url{https://math.berkeley.edu/~mchrist/preprints.html}).

\bibitem{Gr00}
A. Gr\"unrock,
{\it On the Cauchy- and periodic boundary value problem for a certain class of derivative nonlinear Schr\"odinger equations},
arXiv:math/0006195.

\bibitem{Hir14}
H. Hirayama, 
{\it Well-posedness and scattering for a system of quadratic derivative nonlinear Schr\"odinger equations with low regularity initial data}, 
Commun. Pure Appl. Anal. {\bf 13} (2014), no.~4, 1563--1591.

\bibitem{HKO22}
H.~Hirayama, S.~Kinoshita, M.~Okamoto,
\textit{A remark on the well-posedness for a system of quadratic derivative nonlinear Schr\"odinger equations},
Commun. Pure Appl. Anal. \textbf{21} (2022), no. 10, 3309--3334.

\bibitem{HKO25}
H.~Hirayama, S.~Kinoshita, M.~Okamoto,
\textit{Well-posedness and ill-posedness for a system of periodic quadratic derivative nonlinear Schr\"odinger equations},
Pure and Applied Analysis \textbf{7} (2025), no.~2, 359--412.

\bibitem{JeOh26}
I. J. Jeong, S. J. Oh,
{\it Ill-posedness for dispersive equations: degenerate dispersion and the
Takeuchii-Mizohata condition},
Anal. PDE \textbf{19} (2026), no.~4, 721--782.

\bibitem{KO}
T.~Kondo, M.~Okamoto,
\textit{Well- and ill-posedness of the Cauchy problem for semi-linear Schr\"odinger equations on the torus},
to appear in Funkcial. Ekvac.

\bibitem{KoOk26}
T. Kondo, M. Okamoto, 
{\it Norm inflation for quadratic derivative fractional nonlinear Schr\"odinger equations}, 
J. Evol. Equ. {\bf 26} (2026), no.~2, Paper No. 74, 34 pp.

\bibitem{KoOk26b}
T. Kondo, M. Okamoto, 
{\it A sharp criterion for well-posedness of semilinear Schrödinger equations on the real line}, 
arXiv:2610.21734.

\bibitem{Mi85}
S. Mizohata,
{\it On the Cauchy Problem},
Notes and Reports in Mathematics in Science and Engineering, vol.~3,
Academic Press, Orlando, 1985.


\end{thebibliography}
\end{document}